\documentclass[a4paper,12pt,leqno]{amsart}
\usepackage{amsmath,amsthm,amsfonts,amssymb}
\usepackage[titletoc,toc,title]{appendix}
\numberwithin{equation}{section}

\newcommand{\calH}{\mathcal{H}}

\newcommand{\Comp}{\mathbb{C}}

\DeclareMathOperator{\End}{End}
\renewcommand{\epsilon}{\varepsilon}
\DeclareMathOperator{\Id}{Id}
\newcommand{\Integer}{\mathbb{Z}}

\newcommand{\Real}{\mathbb{R}}
\newcommand{\tensor}{\otimes}
\DeclareMathOperator{\tr}{tr}

\newcommand\lemref[1]{Lemma~\ref{#1}}

\newcommand\secref[1]{\S\ref{#1}}
\newcommand\appref[1]{Appendix \ref{#1}}

\newtheorem{theorem}{Theorem}[section]

\newtheorem{lemma}[theorem]{Lemma}

\theoremstyle{definition}

\theoremstyle{remark} 
\begin{document}
\title[$Q$-operators for higher spin eight vertex models]%
{$Q$-operators for higher spin \\eight vertex models\\
with an even number of sites}
\author{Takashi TAKEBE}
\address{Faculty of Mathematics,
National Research University Higher School of Economics,
Vavilova Street 7, Moscow, 117312, Russia.}
\email{\tt ttakebe@hse.ru}

\dedicatory{Dedicated to Professor Evgeny Sklyanin\\
on the occasion of his sixtieth birthday. }
\maketitle
\begin{abstract}
 We construct the $Q$-operator for generalised eight vertex models
 associated to higher spin representations of the Sklyanin algebra,
 following Baxter's 1973 paper. As an application, we prove the sum rule
 for the Bethe roots.
\end{abstract}

\section{Introduction}
\label{sec:intro}

The $Q$-operator was introduced by Baxter \cite{bax:72} in 1972 as an
auxiliary tool to find eigenvalues of the transfer matrix of the eight
vertex model. It satisfies the {\em $TQ$-relation},
\begin{equation*}
    T(u) Q(u) = Q(u) T(u) = h_-(u) Q(u-2\eta) + h_+(u) Q(u+2\eta),
\end{equation*}
with the transfer matrix $T(u)$ and commutes with itself:
$[Q(u),Q(u')]=0$. (The functions $h_\pm(u)$ are defined in
\eqref{def:h+-}.) Under the assumption of semisimplicity of $T(u)$ and
$Q(u)$ the $TQ$-relation leads to equations for the zeros of $Q(u)$,
which eventually give the eigenvalues of $T(u)$. In 1973 Baxter
constructed another $Q$-operator for the eight vertex model with an even
number of sites in \cite{bax:73-1} and explained the relation with the
Bethe Ansatz method in \cite{bax:73-3}. Since then a huge number of
works have been devoted to construction and analysis of $Q$-operators of
various models.

The goal of the present paper is to generalise Baxter's construction of
the $Q$-operator in \cite{bax:73-1} to the higher spin case.

In \cite{tak:92} a generalisation of the eight vertex model was proposed
by means of higher spin representations of the Sklyanin algebra
introduced by Sklyanin in \cite{skl:83}. (The eight vertex model is the
spin $1/2$ case.) It was shown that the algebraic Bethe Ansatz by
Takhtajan and Faddeev \cite{takh-fad:79}, which is an elegant
reformulation of Baxter's construction of the eigenvectors in
\cite{bax:73-1}, \cite{bax:73-2}, \cite{bax:73-3}, can be applied to the
higher spin case as well.

This model was further studied in \cite{tak:95} and \cite{tak:96} by the
algebraic Bethe Ansatz but a certain property of Bethe roots was left
unproved. Baxter showed in \cite{bax:72} that the sum of Bethe roots
satisfies integrality condition, using holomorphicity of the
$Q$-operator. A similar property for the higher spin generalisation was
conjectured in \cite{tak:95} and proved under several ad hoc
assumptions. One of the motivations of the construction of the
$Q$-operator for the generalised models is to prove this sum rule.

\medskip
Remarkably, although Baxter's construction of the $Q$-operator
(\cite{bax:73-1}, \S10.5 of \cite{bax:82}) seems to depend heavily on
the matrix structure of the transfer matrix of the eight vertex model,
it turns out that we can apply his method to our model mutatis
mutandis. In fact,
\begin{itemize}
 \item the auxiliary matrix $M_\lambda(v)$ \eqref{def:gauge-trans} which
       makes an off-diagonal block of the $L$-matrix degenerate is the
       same for any spin and thus the same as that for the eight vertex
       model;
 \item the null vector of the degenerate off-diagonal block of the
       twisted $L$-matrix is the local pseudo vacuum vector (or the
       intertwining vector) $\omega_\lambda(u;v)$ \eqref{pseudo-vac}
       used in the algebraic Bethe Ansatz \cite{tak:92}, \cite{tak:96},
       as is the case with the eight vertex model;
 \item the Hermitian conjugate with respect to the Sklyanin form
       \eqref{skl-form} introduced in \cite{skl:83} plays the role of
       the transpose of matrices in Baxter's construction;
 \item the commutation relation \eqref{QLQR=QLQR} of $Q_R$ and $Q_L$,
       intermediate objects in the construction of the $Q$-operator, is
       proved by using Baxter's argument with the help of Rosengren's
       theories on the elliptic $6j$-symbols \cite{ros:07} and the
       Sklyanin form (the Sklyanin invariant integration) \cite{ros:04}.
\end{itemize}

This paper is organised as follows: in \secref{sec:model} we recall the
model introduced in \cite{tak:92} by defining its transfer matrix. The
construction of the $Q$-operator in \secref{sec:construction} is divided
into four steps as in \cite{bax:82}. As a first step, in
\secref{subsec:QR} we construct an operator $Q_R(u)$ satisfying the half
of the $TQ$-relation by using the local pseudo vacuum vectors. The other
half of the $TQ$-relation (or the ``$QT$''-relation) is satisfied by
another operator $Q_L(u)$, which is defined by the Hermitian conjugate
of $Q_R(u)$ in \secref{subsec:QL}. The most complicated part is the next
step in \secref{subsec:QLQR=QLQR}, the proof of the commutation relation
\eqref{QLQR=QLQR} of $Q_R(u)$ and $Q_L(u)$. The main lemma for this
proof is proved in \appref{app:factorisation}. Once the commutation
relation is proved, the rest of the construction in \secref{subsec:Q} is
a routine work. As an application of the $Q$-operator, we rederive the
Bethe Ansatz equation in \cite{tak:92} and prove the sum rule of the
Bethe roots in \secref{sec:eig-val}. We make several concluding remarks
with comments on related works in the final section
\secref{sec:conclusion}. All the necessary facts about the Sklyanin
algebra are collected in \appref{app:sklyanin}.

\subsection*{Notations}

Throughout this paper we use the following notations and symbols.
\begin{itemize}
 \item $N\in 2\Integer_{>0}$: the number of sites. We consider only even
       $N$. 
 \item $l\in \frac{1}{2}\Integer_{>0}$: the spin of the representation at
       each site. 
 \item $\tau \in i \Real_{>0}$; the elliptic modulus, which is purely
       imaginary. 
 \item $\eta \in [-1/2(2l+1),1/2(2l+1)]$: anisotropy parameter. 
 \item The notations for the theta functions are the same as those in
       Sklyanin's papers \cite{skl:82}, \cite{skl:83} (cf.\
       \cite{mum:83}): 
\begin{equation*}
    \theta_{ab}(z,\tau) = \sum_{n\in\Integer}
       \exp\left( 
             \pi i \left( \frac{a}{2} + n \right)^2 \tau
           +2\pi i \left( \frac{a}{2} + n \right)
                   \left( \frac{b}{2} + z \right)
           \right).
\end{equation*}
       (cf.\ Jacobi's notation (e.g., \cite{whi-wat}): 
       $\vartheta_1(\pi z,\tau)=-\theta_{11}(z,\tau)$, 
       $\vartheta_2(\pi z,\tau)=\theta_{10}(z,\tau)$, 
       $\vartheta_3(\pi z,\tau)=\theta_{00}(z,\tau)$, 
       $\vartheta_4(\pi z,\tau)=\theta_{01}(z,\tau)$.)
 \item We denote $\theta_{11}(z,\tau)$ by $[z]$ for simplicity. 
 \item $[z]_k := \prod_{j=0}^{k-1}[z+2j\eta] =
       [z][z+2\eta]\dotsb[z+2(k-1)\eta]$ for $k=1,2,\dotsc$, $[z]_0=1$.
 \item The Pauli matrices are defined as usual:
\begin{equation*}
    \sigma^0 = \begin{pmatrix} 1 & 0  \\  0  & 1 \end{pmatrix},\quad
    \sigma^1 = \begin{pmatrix} 0 & 1  \\  1  & 0 \end{pmatrix},\quad
    \sigma^2 = \begin{pmatrix} 0 & -i \\  i  & 0 \end{pmatrix},\quad
    \sigma^3 = \begin{pmatrix} 1 & 0  \\  0  &-1 \end{pmatrix}.
\end{equation*}
\end{itemize}

\section{Definition of the model and the $Q$-operator}
\label{sec:model}

In this section we define the generalisation of the eight
vertex model, using the higher spin representations of the Sklyanin
algebra. We mainly follow \cite{tak:95}, with slightly different
normalisation. 

We fix a half integer $l\in\frac{1}{2}\Integer_{>0}$ and consider the
spin $l$ representation space of the Sklyanin algebra as the local state
space $V_i$ ($i=1,\dotsc,N$): $V_i \cong \Theta_{00}^{4l+}$. (See
\appref{app:sklyanin} for the Sklyanin algebra and its representations.)
The total Hilbert space $\calH$ is the tensor product of them:
\begin{equation}
    \calH := V_N \tensor \dotsb \tensor V_1,
\label{def:H}
\end{equation}
and the auxiliary space $V_0$ is a two-dimensional space:
$V_0\cong\Comp^2$. 

The transfer matrix $T(u)$ of the model acting on $\calH$ is defined as
follows: 
\begin{equation}
    T(u) :=
    \tr_0 L_N(u) L_{N-1}(u) \dotsb L_1(u),
\label{def:T(u)}
\end{equation}
where the $L$-operator $L_j(u) \in \End_\Comp(\calH\tensor V_0)$ is
defined as 
\begin{equation}
    L_j(u) = \sum_{a=0}^3 W^L_a(u) \rho_i(S^a) \tensor \sigma^a,
\label{def:Lj}
\end{equation}
by a representation of the Sklyanin algebra on $\calH$,
$
    \rho_i := 1 \tensor \dotsb \tensor
    \rho^l \tensor \dotsb \tensor 1,
$
which acts non-trivially only on $V_i$. The functions
$W^L_a(u)$ in \eqref{def:Lj} are defined by \eqref{def:Wa} in \appref{app:sklyanin} and
the matrices $\sigma^a$ are the Pauli matrices acting on $V_0$.

When the spin $l$ is $1/2$, the local quantum space
$V_i\cong\Theta^{2+}_{00}$ can be identified with $\Comp^2$ by means of
the basis
$
(
    \theta_{00}(2z,2\tau)-\theta_{10}(2z,2\tau),
    \theta_{00}(2z,2\tau)+\theta_{10}(2z,2\tau)
)
$,
and $\rho^{1/2}(S^a)$ are proportional to the Pauli matrices $\sigma^a$,
\eqref{rep:pauli}. Therefore the transfer matrix $T(u)$ is essentially
that of the eight vertex model. The $RLL$-relation \eqref{RLL} leads to
the commutativity of the transfer matrix by the standard argument:
\begin{equation}
    T(u) T(u') = T(u') T(u).
\label{TT}
\end{equation}

The study of the model means the analysis of the spectrum of $T(u)$. The
generalised algebraic Bethe Ansatz (\cite{takh-fad:79}) has been
successfully applied to this model in \cite{tak:92}, \cite{tak:95} and
\cite{tak:96}. However, when Baxter first solved the eight vertex model,
he used the $Q$-operator instead. 

The $Q$-operator is an invertible operator $Q(u):\calH\to\calH$, which
is an entire function in a complex parameter $u$, satisfying the
commutation relations:
\begin{align}
    T(u) Q(u) &= h_-(u) Q(u-2\eta) + h_+(u) Q(u+2\eta),
\label{TQ}
\\
    Q(u) T(u) &= h_-(u) Q(u-2\eta) + h_+(u) Q(u+2\eta),
\label{QT}
\\
   Q(u) Q(v) &= Q(v) Q(u),
\label{QQ=QQ}
\end{align}
where the functions $h_\pm(u)$ are defined by
\begin{equation}
    h_\pm(u) := (2[u\mp2l\eta])^N.
\label{def:h+-}
\end{equation}
Since $T$ is expressed in terms of $Q$, the relation \eqref{QQ=QQ}
implies $[T(u),Q(v)]=0$. We write down equivalent relations \eqref{TQ}
and \eqref{QT} separately, since they arise independently in the
construction.

\section{Construction of the $Q$-operator}
\label{sec:construction}

In this section we construct the $Q$-operator for the higher spin eight
vertex model defined in \secref{sec:model}. The main strategy is the
same as that in Baxter's 1973 paper \cite{bax:73-1}:
\begin{enumerate}
 \item Find an auxiliary two-by-two matrix, a transformation
       by which makes the $(2,1)$-component of the $L$-operator
       degenerate. The tensor products of the null vectors of the
       twisted $(2,1)$-component are the column vectors of the
       $Q_R$-operator satisfying the $TQ$-relation \eqref{TQ}.
 \item Transposing $Q_R$ in an appropriate sense, we obtain the
       $Q_L$-operator which satisfies the $QT$-relation \eqref{QT}.
 \item Show the commutativity $Q_L(u)Q_R(v)=Q_L(v)Q_R(u)$.
 \item The $Q$-operator is defined by
       $Q(u)=Q_R(u)Q_R(u_0)^{-1}=Q_L(u_0)^{-1}Q_L(u)$, where $u_0$ is a
       suitably fixed parameter.
\end{enumerate}

\subsection{Column vectors of $Q_R$}
\label{subsec:QR}

When we applied the generalised algebraic Bethe Ansatz of Takhtajan and
Faddeev \cite{takh-fad:79} in \cite{tak:92}, \cite{tak:95} and
\cite{tak:96}, we used the {\em gauge transformation matrix} of the
$L$-operator of the form\footnote{We use a different normalisation from
those used in \cite{tak:92}, \cite{tak:95}. The matrix $M_\lambda(u)$
here is the same as those used in \cite{tak:97} up to the shift of the
parameter $\lambda\mapsto\lambda+1$ and the multiplication of a diagonal
matrix from the right.}
\begin{equation}
    M_\lambda(v) :=
    \begin{pmatrix}
             - \theta_{00}\left((\lambda - v)/2 , \tau/2\right) &
             - \theta_{00}\left((\lambda + v)/2 , \tau/2\right) \\
    \phantom{-}\theta_{01}\left((\lambda - v)/2 , \tau/2\right) &
    \phantom{-}\theta_{01}\left((\lambda + v)/2 , \tau/2\right)
    \end{pmatrix}.
\label{def:gauge-trans}
\end{equation}
We denote the components of the twisted $L$-operator as
follows\footnote{The spectral parameters $u$ and $v$ here corresponds to
$u-v$ and $u$ in (1.18) of \cite{tak:97}.}.
\begin{equation}
    L_{\lambda,\lambda'}(u;v) =
    \begin{pmatrix}
    \alpha_{\lambda,\lambda'}(u;v) &  \beta_{\lambda,\lambda'}(u;v) \\ 
    \gamma_{\lambda,\lambda'}(u;v) & \delta_{\lambda,\lambda'}(u;v)
    \end{pmatrix}
    := M_\lambda(v)^{-1} L(u) M_{\lambda'}(v),
\label{def:twisted-L}
\end{equation}
where $L(u)$ is the non-trivial part of $L_j(u)$ \eqref{def:Lj}, i.e.,
the operator on $\Theta^{4l+}_{00}\tensor\Comp^2$ defined by
\begin{equation}
    L(u) = \sum_{a=0}^3 W^L_a(u) \rho^l(S^a) \tensor \sigma^a.
\label{def:L}
\end{equation}
It was shown in \cite{tak:92} (see also \cite{tak:95}, \cite{tak:97})
that the operator $\gamma_{\lambda+4l\eta,\lambda}(u;v)$ degenerates and
the following vector $\omega_\lambda(u;v)\in\Theta^{4l+}_{00}$ is a null
vector: 
\begin{equation}
    \omega_\lambda(u;v)
    :=
    \left[
      z + \tfrac{\lambda+u-v}{2}+(-l+1)\eta
    \right]_{2l}
    \left[
     -z + \tfrac{\lambda+u-v}{2}+(-l+1)\eta
    \right]_{2l},
\label{pseudo-vac}
\end{equation}
which is called the {\em local pseudo vacuum} in \cite{tak:97}. This is
a special case ($m=l$) of the {\em intertwining vectors} in
$\Theta^{4l+}_{00}$ defined in \cite{tak:95} (see also \cite{tak:97},
\cite{ros:04} and \cite{ros:07})\footnote{As we changed the
normalisation of $M_\lambda(v)$ in \eqref{def:gauge-trans}, the
normalisation of intertwining vectors is different from that in
\cite{tak:97} by shift $\lambda\mapsto\lambda+1$.}:
\begin{equation}
 \begin{split}
    \phi_{\lambda,\lambda'}(u;v)
    &:=
    \left[
      z + \tfrac{\lambda+u-v}{2}+(-l+1)\eta
    \right]_{l+m}
    \left[
     -z + \tfrac{\lambda+u-v}{2}+(-l+1)\eta
    \right]_{l+m} \times
\\
    &\times
    \left[
      z + \tfrac{\lambda'+u-v}{2}+(-l+1)\eta
    \right]_{l-m}
    \left[
     -z + \tfrac{\lambda'+u-v}{2}+(-l+1)\eta
    \right]_{l-m},
 \end{split}
\label{intertwining-vec}
\end{equation}
where $\lambda'=\lambda+4m\eta$ ($m\in{-l,l+1,\dotsc,l}$). The action
of each component of $L_{\lambda,\lambda'}(u;v)$ on
$\phi_{\lambda',\lambda}(u;v)$ is known (cf.\ \cite{tak:95} (2.4) or
\cite{tak:97} (1.19--22)), among which we need the following three (cf.\
\cite{tak:95} (2.8--10) or 
\cite{tak:97} (1.24--26)):
\begin{equation}
\begin{aligned}
    \alpha_{\lambda+4l\eta,\lambda}(u;v) \omega_{\lambda}(u;v)
    &= 2 [u + 2l\eta] \omega_{\lambda-2\eta}(u;v),
\\
    \gamma_{\lambda+4l\eta,\lambda}(u;v) \omega_{\lambda}(u;v)
    &= 0,
\\
    \delta_{\lambda+4l\eta,\lambda}(u;v) \omega_{\lambda}(u;v)
    &=
    \frac{2[u-2l\eta][\lambda]}{[\lambda+4l\eta]}
    \omega_{\lambda+2\eta}(u;v).
\end{aligned}
\label{action-on-vac}
\end{equation}

An important observation here is that the evenness of $\theta_{00}$ and
$\theta_{01}$ implies the evenness of $M_\lambda(v)$:
\begin{equation}
    M_{-\lambda}(-v) = M_\lambda(v).
\label{M-lambda(-v)}
\end{equation}
Therefore, along with the relations \eqref{action-on-vac}, we have
\begin{equation}
\begin{aligned}
    \alpha_{\lambda-4l\eta,\lambda}(u;v) \omega_{-\lambda}(u;-v)
    &= 2 [u + 2l\eta] \omega_{-\lambda-2\eta}(u;-v),
\\
    \gamma_{\lambda-4l\eta,\lambda}(u;v) \omega_{-\lambda}(u;-v)
    &= 0,
\\
    \delta_{\lambda-4l\eta,\lambda}(u;v) \omega_{-\lambda}(u;-v)
    &=
    \frac{2[u-2l\eta][\lambda]}{[\lambda-4l\eta]}
    \omega_{-\lambda+2\eta}(u;-v).
\end{aligned}
\label{action-on-vac:-}
\end{equation}
Moreover, because of the structure of the local pseudo vacuum vector
\eqref{pseudo-vac}, the shift of the auxiliary parameter $\lambda$ of
$\omega_\lambda(u;v)$ is equivalent to the shift of the spectral
parameter $u$:
\begin{equation}
    \omega_{\lambda\pm2\eta}(u;v) = \omega_{\lambda}(u\pm 2\eta;v).
\label{lamba-shift=u-shift}
\end{equation}
Hence, the relations \eqref{action-on-vac} and \eqref{action-on-vac:-}
mean that, roughly speaking, the operator $\alpha$ (resp.\ $\delta$)
shifts $u$ to $u-2\eta$ (resp.\ $u+2\eta$).

Combining these facts, we construct the vector
$\phi(u;v,\lambda,\vec\sigma)$, which will be a column vector of the
operator $Q_R(u)$. Let us fix complex parameters $v$ and $\lambda$ and a
sequence of $\pm1$,
$\vec\sigma=(\sigma_N,\sigma_{N-1},\dotsc,\sigma_1)$, which satisfies
\begin{equation}
    \sum_{k=1}^N\sigma_k=0.
\label{sum-sigma=0}
\end{equation}
(Since $N$ is even, there are $\binom{N}{N/2}$ sequences satisfying this
condition.)  We define the sequence $\{\lambda_j\}_{j=1,\dots,N+1}$ by
\begin{equation}
    \lambda_1 = \lambda, \quad
    \lambda_{j+1} := \lambda_j + 4\sigma_j l \eta
              = \lambda + 4l\eta \sum_{k=1}^{j}\sigma_k.
\label{def:lambdaj}
\end{equation}
The condition \eqref{sum-sigma=0} implies
$\lambda_{N+1}=\lambda_1$. Because of the evenness \eqref{M-lambda(-v)}
of $M_\lambda(v)$, we have
\begin{equation}
    M_{\sigma_j \lambda_j+4l\eta}(\sigma_j v)
    =
    M_{\lambda_j+ 4 \sigma_j l\eta}(v)
    =
    M_{\lambda_{j+1}}(v)
    =
    M_{\sigma_{j+1} \lambda_{j+1}}(\sigma_{j+1} v).
\label{M(j+1)=M(j,4leta)}
\end{equation}
($\sigma_{N+1}:=\sigma_1$.)

Therefore, if we define a vector $g_j(u)=g_j(u;v,\lambda,\vec\sigma)$ in
$V_j$ by 
\begin{equation}
    g_j(u;v,\lambda,\vec\sigma)
    :=
    \omega_{\sigma_j\lambda_j}(u;\sigma_j v),
\label{def:gj}
\end{equation}
then the formulae \eqref{action-on-vac} and \eqref{action-on-vac:-}
imply  
\begin{equation}
 \begin{split}
    \alpha_{\lambda_{j+1},\lambda_j}(u;v) g_j(u)
    &=
    \alpha_{\lambda_j+4\sigma_j l\eta,\lambda_j}(u;v)
    \omega_{\sigma_j\lambda_j}(u;\sigma_j v)
\\
    &= 2 [u + 2l\eta] g_j(u-2\eta),
\\
    \gamma_{\lambda_{j+1},\lambda_j}(u;v) g_j(u)
    &=
    \gamma_{\lambda_j+4\sigma_j l\eta,\lambda_j}(u;v)
    \omega_{\sigma_j \lambda}(u;\sigma_j v)
\\
    &= 0,
\\
    \delta_{\lambda_{j+1},\lambda_j}(u;v) g_j(u)
    &=
    \delta_{\lambda_j+4l\eta,\lambda_j}(u;v)
    \omega_{\sigma_j \lambda_j}(u;\sigma_j v)
\\
    &=
    \frac{2[u-2l\eta][\lambda_j]}{[\lambda_{j+1}]}
    g_j(u+2\eta).
 \end{split}
\label{action-on-g}
\end{equation}

Let us compute the action of the transfer matrix on the tensor product 
\begin{equation}
    \phi(u;v,\lambda,\vec\sigma)
    :=
    g_N(u;v,\lambda,\vec\sigma) \tensor 
    g_{N-1}(u;v,\lambda,\vec\sigma) \tensor \dotsb \tensor
    g_1(u;v,\lambda,\vec\sigma) \in\calH.
\label{def:phi}
\end{equation}
By insertion of
$
    1
    = 
    M_{\lambda_j}(v) M_{\lambda_j}(v)^{-1}
$
between $L_j(u)$ and $L_{j-1}(u)$ in the definition \eqref{def:T(u)} of
the transfer matrix and by the cyclicity of the trace, we can rewrite
the transfer matrix as
\begin{equation}
    T(u)
    =
    \tr_0
    \prod_{j=1,\dots,N}^{\curvearrowleft}
    \begin{pmatrix}
    \alpha_{\lambda_{j+1},\lambda_j}(u;v) &
    \beta_{\lambda_{j+1},\lambda_j}(u;v) \\
    \gamma_{\lambda_{j+1},\lambda_j}(u;v) &
    \delta_{\lambda_{j+1},\lambda_j}(u;v) 
    \end{pmatrix}
\label{T(u)=prod(twisted-L)}
\end{equation}
The formulae \eqref{action-on-g} reduces
$T(u)\phi(u;v,\lambda,\vec\sigma)$ to a triangular form. Thus we obtain
\begin{equation}
    T(u) \phi(u;v,\lambda,\vec\sigma)
    =
    h_-(u) \phi(u-2\eta;v,\lambda,\vec\sigma)
    +
    h_+(u) \phi(u+2\eta;v,\lambda,\vec\sigma).
\label{T(u)phi}
\end{equation}
Let
$
    \{
    \phi_k(u)
    :=
    \phi(u;v_k,\lambda_k,\vec\sigma_k)
    \}_{k=1,\dotsc,\dim\calH}
$
be a set of $\dim\calH = (2l+1)^N$ vectors with distinct values of
parameters. We define a linear operator
$Q_R(u):\Comp^{\dim\calH}\to\calH$ by
\begin{equation}
    Q_R(u): e_k \mapsto \phi_k(u),
\label{def:QR}
\end{equation}
where $\{e_k\}_{k=1,\dotsc,\dim\calH}$ is a basis of
$\Comp^{\dim\calH}$. The $TQ$-relation \eqref{TQ} for $Q_R(u)$,
\begin{equation}
     T(u) Q_R(u) = h_-(u) Q_R(u-2\eta) + h_+(u) Q_R(u+2\eta),
\label{TQR}
\end{equation}
is a direct consequence of \eqref{T(u)phi}.

\subsection{Hermitian conjugate and $Q_L$}
\label{subsec:QL}

As the next step, we construct an operator $Q_L(u)$ satisfying the
``$QT$''-relation \eqref{QT}. For this purpose we study the Hermitian
adjoint of $T(u)$. The space $\calH$ has a natural Hermitian structure
induced from the Sklyanin form \eqref{skl-form} on
$\Theta^{4l+}_{00}$. We consider the adjoint operator with respect to
this Hermitian structure.

Denote the four component of the $L$-operator \eqref{def:L} by
$L_{\epsilon\epsilon'}(u)$ ($\epsilon, \epsilon' = \pm$):
\begin{equation}
    L(u) =
    \begin{pmatrix}
     L_{--}(u) & L_{-+}(u) \\ L_{+-}(u) & L_{++}(u)
    \end{pmatrix}.
\label{L+-}
\end{equation}
By the definition \eqref{def:Wa} of the functions $W_a^L(u)$ and the
self-adjointness \eqref{Sa:self-adj} of $S^a$, it is easy to see that
the adjoint operator of $L_{\epsilon\epsilon'}$ is expressed by
$L_{-\epsilon,-\epsilon'}$ as follows:
\begin{equation}
 \begin{aligned}
    (L_{--}(u))^* &= - L_{++}(-\bar u), &
    (L_{-+}(u))^* &= \phantom{-} L_{+-}(-\bar u),
\\
    (L_{+-}(u))^* &= \phantom{-} L_{-+}(-\bar u), &
    (L_{++}(u))^* &= - L_{--}(-\bar u).
 \end{aligned}
\label{(L+-)*}
\end{equation}
Substituting the expression \eqref{L+-} into the definition
\eqref{def:T(u)} of the transfer matrix $T(u)$, we have
\begin{equation}
 \begin{split}
    T(u) &=
    \sum_{\epsilon_{N-1},\dotsc,\epsilon_1=\pm}
    L_{N,-\epsilon_{N-1}}(u) \dotsb
    L_{j,\epsilon_{j}\epsilon_{j-1}}\dotsb
    L_{1,\epsilon_{1}-}
\\
    &+
    \sum_{\epsilon_{N-1},\dotsc,\epsilon_1=\pm}
    L_{N,+\epsilon_{N-1}}(u) \dotsb
    L_{j,\epsilon_{j}\epsilon_{j-1}}(u)\dotsb
    L_{1,\epsilon_{1}+}(u),
\end{split}
\label{T=L+L}
\end{equation}
where $L_{j,\epsilon\epsilon'}(u)$ is the component of $L_j(u)$ defined
by \eqref{def:Lj}. Note that each term in the first sum corresponds
bijectively to a term in the second sum by flipping all signs
$\epsilon_j$. On the other hand, there are as many $L_{-+}$ as $L_{+-}$
in each sum because of the boundary condition $\epsilon_N=\epsilon_0$
($=$ the $\pm$ signs of the rightmost and the leftmost factors), which
comes from the trace in \eqref{def:T(u)}. Therefore, in each summand in
\eqref{T=L+L}, 
\begin{multline}
    (\text{number of }L_{j,--}) +
    (\text{number of }L_{j,++})
\\
    =
    N - 2(\text{number of }L_{j,+-})
    \equiv N \equiv 0 \pmod{2}.
\label{number(L--)+number(L++)}
\end{multline}
The adjoint of $T(u)$ is obtained by substitution of \eqref{(L+-)*} into
\eqref{T=L+L}:
\begin{equation}
 \begin{split}
    \bigl(T(u)\bigr)^* &=
    \sum_{\epsilon_{N-1},\dotsc,\epsilon_1=\mp}
    L_{N,+\epsilon_{N-1}}(-\bar u) \dotsb
    L_{j,\epsilon_{j}\epsilon_{j-1}}(-\bar u)\dotsb
    L_{1,\epsilon_{1}+}(-\bar u)
\\
    &+
    \sum_{\epsilon_{N-1},\dotsc,\epsilon_1=\mp}
    L_{N,-\epsilon_{N-1}}(-\bar u) \dotsb
    L_{j,\epsilon_{j}\epsilon_{j-1}}(-\bar u)\dotsb
    L_{1,\epsilon_{1}-}(-\bar u)
\\
    &= T(-\bar u).
\end{split}
\label{T*}
\end{equation}
The signs in \eqref{(L+-)*} vanish in \eqref{T*} because of
\eqref{number(L--)+number(L++)}. 

Now let us take the adjoint of \eqref{TQR} with $-\bar u$ instead of $u$. We
have
\begin{equation*}
 \begin{split}
    Q_R(-\bar u)^* T(u)
    &=
    \overline{h_-(-\bar u)} Q_R(-\bar u - 2\eta)^*
    +
    \overline{h_+(-\bar u)} Q_R(-\bar u + 2\eta)^*
\\
    &=
    h_+(u) Q_R(-\bar u - 2\eta)^*
    +
    h_-(u) Q_R(-\bar u + 2\eta)^*,
 \end{split}
\end{equation*}
because the theta function $[u]=\theta_{11}(u,\tau)$ is an odd real
analytic function. (Here again we use the fact that $N$ is even.)
Defining an operator $Q_L(u)$ by
\begin{equation}
    Q_L(u) := Q_R(-\bar u)^*:\calH \to \Comp^{\dim\calH},
\label{def:QL}
\end{equation}
we obtain the operator satisfying the relation \eqref{QT}:
\begin{equation}
    Q_L(u) T(u) = h_-(u) Q_L(u-2\eta) + h_+(u) Q_L(u+2\eta).
\label{QLT}
\end{equation}

\subsection{Commutation relation of $Q_R$ and $Q_L$}
\label{subsec:QLQR=QLQR}

An important property of the operators $Q_R(u)$ and $Q_L(u)$
constructed above, \eqref{def:QR} and \eqref{def:QL}, is the commutation
relation: 
\begin{equation}
    Q_L(u) Q_R(u') = Q_L(u') Q_R(u).
\label{QLQR=QLQR}
\end{equation}
The $(i,j)$-element of $Q_L(u) Q_R(u')$ is
\begin{equation*}
 \begin{split}
    (e_i,Q_L(u) Q_R(u')e_j) &=
    (e_i, Q_R(-\bar u)^* Q_R(u') e_j)
\\
    &=
    \langle Q_R(-\bar u) e_i,  Q_R(u') e_j\rangle
    = \langle \phi_i(-\bar u), \phi_j(u')\rangle.
 \end{split}
\end{equation*}
where $(,)$ is an Hermitian form on $\Comp^{\dim\calH}$ defined by
$(e_i,e_j)=\delta_{ij}$. Hence, if the function $\Phi(u,u')$ of $(u,u')$
defined by
\begin{equation}
    \Phi(u,u')
    :=
    \langle
     \phi(-\bar u; v,\lambda,\vec\sigma),
     \phi(     u'; v',\lambda',\vec\sigma')
    \rangle
\label{def:Phi}
\end{equation}
is symmetric in $u$ and $u'$ for any choice of parameters
$(v,\lambda,\vec\sigma)$ and $(v',\lambda',\vec\sigma')$, the
commutation relation \eqref{QLQR=QLQR} holds. By the definitions
\eqref{def:phi} and \eqref{def:gj}, $\Phi(u,u')$ is rewritten as
\begin{equation}
    \Phi(u,u')
    =
    \prod_{k=1}^N
    \langle
     \omega_{\sigma_k  \lambda_k }(-\bar u;\sigma_k  v ),
     \omega_{\sigma_k' \lambda_k'}(     u';\sigma_k' v')
    \rangle.
\label{Phi=prod}
\end{equation}

From the results of Rosengren \cite{ros:04} and \cite{ros:07} follows
the factorisation of $\Phi(u,u')$.
\begin{lemma}
\label{lem:factorise-Phi}
 $\Phi(u,u')$ is factorised as
\begin{equation}
 \begin{split}
    \Phi(u,u') ={}& C \prod_{k=1}^N
    F\bigl(
      \tfrac{\lambda'_k - \bar\lambda_k}{2}
     +\tfrac{\sigma'_k u' + \sigma_k u}{2}
     +(\sigma'_k-\sigma_k)l\eta
     +\tfrac{-v'+\bar v}{2}
     \bigr) \times
\\
    &\times \prod_{k=1}^N
    G\bigl(
      \tfrac{\lambda'_k + \bar\lambda_k}{2}
     +\tfrac{\sigma'_k u' - \sigma_k u}{2}
     +(\sigma'_k-\sigma_k)l\eta
     +\tfrac{-v'-\bar v}{2}
     \bigr).
 \end{split}
\label{Phi=FG}
\end{equation}
 Here the constant $C$ depends only on $\tau$, $\eta$, $l$ and $N$,
 while the functions $F(w)$ and $G(w)$ depend only on $\tau$, $\eta$ and
 $l$.
\end{lemma}
We prove this lemma in \appref{app:factorisation}.

This factorisation is exactly of the same type as (10.5.27) in
\cite{bax:82} for the eight vertex model. The parameters $\lambda$ and
$s_j$ in \cite{bax:82} correspond to $2l\eta$ and $\lambda_j/2$ here
respectively. Baxter showed (p.219 (\S10.5) of \cite{bax:82}) that a
function with such a factorisation is symmetric in $u$ and $u'$,
whatever $F$ and $G$ are.

Thus the commutation relation \eqref{QLQR=QLQR} has been proved.

\subsection{The $Q$-operator and its commutation relations}
\label{subsec:Q}

The rest of the construction of the $Q$-operator is the same as that in
\cite{bax:73-1} and \cite{bax:82}. 

As in \cite{bax:73-1} and \cite{bax:82}, we expect that varying parameters
$(v,\lambda,\vec\sigma)=(v_k,\lambda_k,\vec\sigma_k)$ produces a set of
sufficiently many vectors
$
    \{
    \phi_k(u)
    :=
    \phi(u;v_k,\lambda_k,\vec\sigma_k)
    \}_{k=1,\dotsc,\dim\calH}
$
which spans the space $\calH$ for generic $u$, although we do not have a
proof. Here we suppose this non-degeneracy and take a parameter $u=u_0$,
for which $Q(u_0)$ has the inverse. We define the $Q$-operator by
\begin{equation}
    Q(u) := Q_R(u) Q_R(u_0)^{-1} : 
    \calH \xrightarrow{Q_R(u_0)^{-1}}
    \Comp^{\dim\calH} \xrightarrow{Q_R(u)} \calH. 
\label{Q=QRQR}
\end{equation}

Then the commutation relation \eqref{QLQR=QLQR} implies
\begin{gather}
    Q(u) = Q_L(u_0)^{-1} Q_L(u),
\label{Q=QLQL}
\\
    Q(u) Q(u')
    = Q_L(u_0)^{-1} Q_L(u) Q_R(u') Q_R(u_0)^{-1}
    = Q(u') Q(u). 
\label{QQ=QQ:proof}
\end{gather}
Multiplying $Q_R(u_0)^{-1}$ from the right to \eqref{TQR}, we have
\begin{equation*}
    T(u) Q(u) = h_-(u) Q(u-2\eta) + h_+(u) Q(u+2\eta),
\end{equation*}
while left multiplication of $Q_L(u_0)^{-1}$ to \eqref{QLT} leads to 
\begin{equation*}
    Q(u) T(u) = h_-(u) Q(u-2\eta) + h_+(u) Q(u+2\eta).
\end{equation*}
Thus we have all the commutation relations \eqref{TQ}, \eqref{QT},
\eqref{QQ=QQ} and
\begin{equation}
    T(u)Q(u)=Q(u)T(u)
\label{TQ=QT}
\end{equation}
from \eqref{TQ} and \eqref{QT}. 

\section{The eigenvalue of the transfer matrix}
\label{sec:eig-val}

As an application of the $Q$-operator, let us compute the eigenvalue of
the transfer matrix $T(u)$ and prove the sum rule of the Bethe roots,
which was proved in \cite{tak:95} under redundant conditions.

For the eight vertex model there are two involutive operators on $\calH$
($R =$ (reversing the arrows) and $S =$ (assigning $-1$ to down arrows)
in \cite{bax:82}) which commute with the transfer matrix. Baxter used
them to break up $\calH$ and deduced the Bethe Ansatz equation.

In our case we have also two involutions, $U_1^{\tensor N}$ and
$U_3^{\tensor N}$, corresponding to $R$ and $S$. The definition of
$U_a$, which was introduced by Sklyanin in \cite{skl:83}, is given by
\eqref{def:Ua} in \appref{app:sklyanin}. It was shown in \S1.4 of
\cite{tak:95} that they commute with each other and with the transfer
matrix $T(u)$.
\begin{equation}
    (U_a^{\tensor N})^2 = 1, \quad
    [U_a^{\tensor N}, U_b^{\tensor N}] = [T(u), U_a^{\tensor N}]=0.
\label{Ua:comm-rel}
\end{equation}

By the definition \eqref{def:Ua} of $U_a$, the operators $U_1$ and $U_3$
act on the local pseudo vacuum vector $\omega_\lambda(u;v)$ defined by
\eqref{pseudo-vac} as
\begin{align}
    U_1 \omega_\lambda(u;v)
    &=
    e^{-l\pi i} \omega_\lambda(u+1;v),
\label{U1-on-omega}
\\
    U_3 \omega_\lambda(u;v)
    &=
    e^{l\pi i (\tau-1) + 2l \pi i (\lambda + u - v + 2l\eta)}
    \omega_\lambda(u+\tau;v).
\label{U3-on-omega}
\end{align}
Therefore operators $U_a^{\tensor N}$ act on the column vector
$\phi(u;v,\lambda,\vec\sigma)$ of the $Q_R$-operator defined by
\eqref{def:phi} as follows:
\begin{equation}
 \begin{aligned}
    U_1^{\tensor N} \phi(u;v,\lambda,\vec\sigma)
    &=
    e^{-N l \pi i} \phi(u+1;v,\lambda,\vec\sigma),
\\
    U_3^{\tensor N} \phi(u;v,\lambda,\vec\sigma)
    &=
    e^{N l \pi i(\tau-1) + 2Nl \pi i u}
    \phi(u+\tau;v,\lambda,\vec\sigma).
 \end{aligned}
\label{Ua-on-phi}
\end{equation}
Here we used a formula (the same as (10.5.40) in \cite{bax:82})
\begin{equation*}
    \sum_{k=1}^N \sigma_k \lambda_k
    =
    -2Nl\eta,
\end{equation*}
which is a consequence of the condition \eqref{sum-sigma=0} and the
definition \eqref{def:lambdaj} of $\lambda_j$. Note that the
coefficients of $\phi$ and the shifts of $u$ in \eqref{Ua-on-phi} do not
depend on the parameters $(v,\lambda,\vec\sigma)$. Therefore the
$Q_R$-operator defined by \eqref{def:QR} inherits those relations:
\begin{align}
    U_1^{\tensor N} Q_R(u)
    &=
    e^{-N l \pi i} Q_R(u+1),
\label{U1-on-QR}
\\
    U_3^{\tensor N} Q_R(u)
    &=
    e^{N l \pi i(\tau-1) + 2Nl \pi i u} Q_R(u+\tau).
\label{U3-on-QR}
 \end{align}
Because of the unitarity of the involution $U_a^{\tensor N}$ and the
definition \eqref{def:QL} of $Q_L(u)$, we have
\begin{equation*}
    Q_L(u)U_a^{\tensor N}
    =
    \bigl( U_a^{\tensor N} Q_R(-\bar u) \bigr)^*.
\end{equation*}
When $a=1$, the product $U_1^{\tensor N} Q_R(-\bar u)$ is equal to
$e^{-Nl\pi i}Q_R(-\bar u+1)$ by \eqref{U1-on-QR}. Note that the property
of the theta function $[z+1]=-[z]$ implies
$\omega_\lambda(u+2;v)=\omega_\lambda(u;v)$ and $Q_R(u+2)=Q_R(u)$. Hence
\begin{equation*}
    U_1^{\tensor N} Q_R(\bar u)
    =
    e^{-Nl\pi i}Q_R(-\bar u -1)
    =
    e^{-Nl\pi i} Q_R(-\overline{(u+1)}),
\end{equation*}
which gives
\begin{equation}
    Q_L(u) U_1^{\tensor N} = e^{-Nl\pi i} Q_L(u+1).
\label{U1-on-QL}
\end{equation}
(Recall that $Nl$ is an integer. Therefore $e^{-Nl\pi i}$ is a real
number $\pm1$.) Likewise we can prove
\begin{equation}
    Q_L(u) U_3^{\tensor N}
    =
    e^{N l \pi i(\tau-1) + 2Nl \pi i u} Q_L(u+\tau).
\label{U3-on-QL}
\end{equation}
Here we used $-\bar u+\tau = -\overline{(u+\tau)}$ and $e^{N l \pi
i(\tau-1)}\in\Real$, which are consequences of $\tau\in i\Real$.

Multiplying $Q_R(u_0)^{-1}$ to \eqref{U1-on-QR} and \eqref{U3-on-QR}
from the right and \ $Q_L(u_0)^{-1}$ to \eqref{U1-on-QL} and
\eqref{U3-on-QL} from the left, we obtain 
\begin{equation}
 \begin{aligned}
    U_1^{\tensor N} Q(u) = Q(u) U_1^{\tensor N}
    &=
    e^{-Nl\pi i} Q(u+1),
\\
    U_3^{\tensor N} Q(u) = Q(u) U_3^{\tensor N}
    &=
    e^{N l \pi i(\tau-1) + 2Nl \pi i u} Q(u+\tau).
 \end{aligned}
\label{Ua-on-Q}
\end{equation}

\bigskip
Having shown that the operators $T(u)$, $Q(u)$ and $U_a^{\tensor N}$
commute with each other, we can consider the diagonalisation problem of
$T(u)$ and $Q(u)$ on the joint eigenspaces of $U_1^{\tensor N}$ and
$U_3^{\tensor N}$. Recall that $(U_a^{\tensor N})^2=1$
\eqref{Ua:comm-rel}. Hence the operators $U_a^{\tensor N}$ are
diagonalisable and their eigenvalues are $\pm1$. Accordingly the space
$\calH$ is decomposed as
\begin{equation}
    \calH = \bigoplus_{\nu_1,\nu_3=0,1} \calH_{\nu_1,\nu_3},
    \qquad
    U_a^{\tensor N}|_{\calH_{\nu_1,\nu_3}}
    = (-1)^{\nu_a} \Id_{\calH_{\nu_1,\nu_3}}.
\label{H=Hnu}
\end{equation}
The rest is the same as in \S10.6 of \cite{bax:82}. We suppose that
$T(u)$ and $Q(u)$ are diagonalisable. Because of the commutativity, they
are diagonalised simultaneously on each $\calH_{\nu_1,\nu_3}$ by a
matrix independent of $u$. Let us denote the eigenvalues of $T(u)$ and
$Q(u)$ for one of the eigenvectors by $\Lambda(u)$ and $q(u)$, which are
entire functions in $u$. The $TQ$-relation \eqref{TQ} gives the relation
\begin{equation}
    \Lambda(u) q(u) = h_-(u) q(u-2\eta) + h_+(u) q(u+2\eta).
\label{tq}
\end{equation}
The transformation rules \eqref{Ua-on-Q} reduce to
the scalar equations:
\begin{equation}
\begin{aligned}
    (-1)^{\nu_1} q(u) &= e^{-Nl\pi i} q(u+1),
\\
    (-1)^{\nu_3} q(u) &= e^{N l \pi i(\tau-1) + 2Nl \pi i u} q(u+\tau).
\end{aligned}
\label{nuq}
\end{equation}
Applying the argument principle in the complex analysis, we can prove
from \eqref{nuq} that the entire function $q(u)$ has $Nl$ zeros in the
rectangular with vertices $0$, $1$, $1+\tau$ and $\tau$. Let us denote
these zeros by $u_1,\dotsc,u_{Nl}$. Note that $q(u)$ has zeros at
$u_j+n+m\tau$ ($n,m\in\Integer$) as well because of the
quasi-periodicity \eqref{nuq}.

The theta function $[z]=\theta_{11}(z,\tau)$ has the quasi-periodicity:
\begin{equation}
    [z+1] = - [z], \qquad
    [z+\tau] = - e^{-\pi i \tau - 2\pi i z}[z],
\label{theta11-period}
\end{equation}
and $[z]=0$ is equivalent to $z\in\Integer+\tau\Integer$. Hence the
function
\begin{equation}
    f(u) := \frac{q(u)}{\prod_{j=1}^{Nl} [u-u_j]}
\label{def:q/theta}
\end{equation}
is entire, does not vanish and has the quasi-periodicity:
\begin{equation*}
    f(u+1) = (-1)^{\nu_1} f(u), \qquad
    f(u+\tau)
    = 
    (-1)^{\nu_3} e^{ - 2\pi i \sum_{j=1}^{Nl}u_j} f(u)
\end{equation*}
because of \eqref{nuq}. By the standard argument in the complex analysis
again, it follows from this periodicity that $f(u)$ is a constant
multiple of $e^{K u}$, where the constant $K$ satisfies
\begin{equation}
    K
    =
    \nu_1 \pi i + 2n_1\pi i, \qquad
    K\tau
    =
    \nu_3 \pi i - 2\pi i \sum_{j=1}^{Nl}u_j + 2n_3\pi i
\label{exp(sum-rule)}
\end{equation}
for some integers $n_1$ and $n_3$. Putting the condition
\eqref{exp(sum-rule)} back into the definition \eqref{def:q/theta} of
$f(u)$, we obtain the explicit form of $q(u)$
\begin{equation}
    q(u)
    = 
    C e^{\nu_1 \pi i u} \prod_{j=1}^{Nl} [u-u_j],
\label{q(u)}    
\end{equation}
and the {\em sum rule}\footnote{The conjecture in \cite{tak:96} should
be modified.}:
\begin{equation}
    \sum_{j=1}^{Nl}u_j 
    \equiv
    -\frac{\nu_1\tau}{2} + \frac{\nu_3}{2}
    \pmod{\Integer + \tau\Integer}.
\label{sum-rule}
\end{equation}
Setting $u=u_j$ in \eqref{tq}, we have the equation
\begin{equation*}
    h_-(u_j) q(u_j-2\eta) + h_+(u_j) q(u_j+2\eta)=0,
\end{equation*}
or equivalently,
\begin{equation}
    \left(\frac{[u_j+2l\eta]}{[u_j-2l\eta]}\right)^N
    =
    e^{4\nu_1 \pi i \eta}
    \prod_{k=1,k\neq j}^{Nl}
    \frac{[u_j-u_k+2\eta]}{[u_j-u_k-2\eta]},
\label{bethe-eq}
\end{equation}
which is nothing but the {\em Bethe equation}. The corresponding
eigenvalue of the transfer matrix is
\begin{equation}
 \begin{split}
    \Lambda(u)
    ={}&
    (2[u+2l\eta])^N e^{-2\nu_1\pi i\eta}
    \prod_{j=1}^{Nl}
    \frac{[u-u_j-2\eta]}{[u-u_j]}
\\
    +{}&
    (2[u-2l\eta])^N e^{ 2\nu_1\pi i\eta}
    \prod_{j=1}^{Nl}
    \frac{[u-u_j+2\eta]}{[u-u_j]},
 \end{split}
\label{bethe-eigen-val}
\end{equation}
which was first obtained in \cite{tak:92} by the Bethe Ansatz.

\section{Concluding comments and remarks}
\label{sec:conclusion}

We have constructed the $Q$-operator (\eqref{Q=QRQR} or \eqref{Q=QLQL})
satisfying the $TQ$-relation, \eqref{TQ} and \eqref{QT}, and commuting
with itself, \eqref{QQ=QQ}, for generalised eight vertex models defined
by the higher spin representations of the Sklyanin algebra. Our explicit
construction by means of the intertwining vectors makes it possible to
prove the sum rule \eqref{sum-rule} of the Bethe roots, which we could
not prove in \cite{tak:95} from the Bethe Ansatz itself. Thus the
$Q$-operator can be useful to analyse the Bethe roots.

However there are several weak points in our construction. Let us make
several comments on them with remarks on recent developments.

\begin{itemize}
 \item Our construction as well as Baxter's original idea relies on the
       assumption that the operator $Q_R(u)$ is generically
       non-degenerate. (See \eqref{Q=QRQR}.) This is very difficult to
       prove, although is plausible, as Baxter claimed on the basis of
       degenerate (six vertex) cases (p.220, \cite{bax:82}).

 \item The Bethe Ansatz for the higher spin generalisation of the eight
       vertex model in \cite{tak:92}, \cite{tak:95} and \cite{tak:96}
       works when the total spin $Nl$ is an integer. In particular, if
       $l$ is an integer, there is no constraint imposed on $N$. On the
       other hand, the construction of the $Q$-operator in this paper
       requires that $N$ is always even, especially because of the
       condition \eqref{sum-sigma=0}, $\sum \sigma_k = 0$. In this
       respect our $Q$-operator method is weaker than the Bethe Ansatz.

 \item As is mentioned above, when the number of the lattice sites is
       odd, our construction does not work. In this case, we would need
       higher spin generalisation of Baxter's 1972 paper
       \cite{bax:72}. There are many works (for example
       \cite{fab-mcc:03-05}, \cite{fab:07}, \cite{fab-mcc:07},
       \cite{roa:07-2}) in this direction for the eight vertex model.

 \item The $Q$-operators for the XXZ spin chain of higher spin were
       constructed by Roan in \cite{roa:07-1}, following both ways of
       Baxter, \cite{bax:72} and \cite{bax:73-1}. The latter method is
       similar to ours, but since all the ingredients in the
       trigonometric case have explicit matrix description, there are no
       complication caused by the functional realisation of the
       representation spaces as in our case. (See also \cite{mot:13} for
       a related work.)

 \item In addition to the difficulty of the non-degeneracy problem of
       $Q_R(u)$, there is another weakness of Baxter's classical method:
       the construction is essentially non-local because of $Q_R^{-1}$
       in \eqref{Q=QRQR}. Therefore it is not very useful for a further
       analysis of the model except its spectrum. Modern developments
       from different viewpoints\footnote{The author thanks the referee
       for informing the references in this remark and \cite{c-d-s:14}
       in the next remark.} remove this weak point. In the context of
       the conformal field theory, Bazhanov, Lukyanov and Zamolodchikov
       (\cite{b-l-z:97}, \cite{b-l-z:99}) suggested an essentially new
       approach, which gives $Q$-operators as traces of monodromy
       matrices in the trigonometric spin$\,=1/2$ case. Mangazeev
       extended this construction to higher spin trigonometric cases in
       \cite{man:14-1} (using the fusion procedure of the transfer
       matrices) and in \cite{man:14-2} (as integral operators), based
       on the ideas of factorised $L$-operators
       \cite{bazh-str:90}. Higher rank generalisation is found in
       \cite{b-f-l-m-s:11}. These methods provide local construciton of
       $Q$-operators.

 \item The $Q$-operators for the elliptic models with
       infinite-dimensional state spaces were constructed by Zabrodin in
       \cite{zab:00} and by Chicherin, Derkachov, Karakhanyan and
       Kirschner in \cite{c-d-k-k:13}. 

       It is an interesting question, whether our $Q$-operator can be
       obtained by reduction from their $Q$-operators. In fact, recently
       Chicherin, Derkachov and Spiridonov \cite{c-d-s:14} constructed a
       general elliptic $R$-operator acting in the tensor product of
       infinite dimensional representations of the Sklyanin algebra as
       the integral operator similar to \cite{man:14-2} and restricted
       it to finite-dimensional representations. It would be remarkable,
       if such a restriction of elliptic $Q$-operators in \cite{zab:00}
       or \cite{c-d-k-k:13} would be found, as they are constructed for
       any (i.e., even and odd) number of sites.

 \item Recently the $Q$-operator is studied from the viewpoint of
       representation theory of quantum algebras. See, for example,
       \cite{fre-her:15}. It is a challenging problem to understand
       elliptic $Q$-operators from the representation theory of elliptic
       quantum algebras.

\end{itemize}

\subsection*{Acknowledgements}

The author dedicates this paper to Professor Evgeny Sklyanin on the
occasion of his sixtieth birthday. It was Professor Sklyanin, who first
led the author to the study of the elliptic quantum integrable models.

\medskip
Special thanks are due to Hitoshi Konno, who suggested to use
Rosengren's results in the proof of \lemref{lem:factorise-Phi}. The
author also expresses his gratitude to Masahiko Ito, Saburo Kakei, Atsuo
Kuniba, Masatoshi Noumi, Jun'ichi Shiraishi and Anton Zabrodin for
discussions and encouragement.

The author is grateful to Rikkyo University, Tokyo, for its
hospitality, where a part of this work was done.

The article was prepared within the framework of the Academic fund
Program at the National Research University Higher School of Economics
(HSE) in 2015--2016 (grant No.15-01-102) and supported within the
framework of a subsidy granted to the HSE by the Government of the
Russian Federation for the implementation of the Global Competitiveness
Program. 

\begin{appendices}
\section{Sklyanin algebra}
\label{app:sklyanin}

In this appendix we recall several facts on the Sklyanin algebra and
its representations from \cite{skl:82} and \cite{skl:83}.

The {\em Sklyanin algebra} is an associative algebra generated by four
generators $S^a$ ($a=0,\dotsc,3$) subject to the following relations:
\begin{equation}
    L_{12}(v) L_{13}(u) R_{23}(u-v) = R_{23}(u-v) L_{13}(u) L_{12}(v).
\label{RLL}
\end{equation}
Here the symbols are defined as follows:
\begin{itemize}
 \item the {\em $L$-operator} $L(u)$ with a complex parameter $u$ is
defined by
\begin{equation}
    L(u) = \sum_{a=0}^3 W_a^L(u) S^a \tensor \sigma^a,
\label{def:L:alg}
\end{equation}
where\footnote{The functions $W_a^L(u)$ here are normalised differently
       from those in \cite{tak:97}: $W_a^{L \text{(here)}}(u) =
       2\theta_{11}(2\eta,\tau)W_ a^{L \text{(old)}}$.}
\begin{equation}
\begin{aligned}
    W_0^L(u) 
    &= \frac{\theta_{11}(u,\tau)}{\theta_{11}(\eta,\tau)},&
    W_1^L(u)
    &= \frac{\theta_{10}(u,\tau)}{\theta_{10}(\eta,\tau)},
\\
    W_2^L(u)
    &= \frac{\theta_{00}(u,\tau)}{\theta_{00}(\eta,\tau)},&
    W_3^L(u)
    &= \frac{\theta_{01}(u,\tau)}{\theta_{01}(\eta,\tau)}.
\end{aligned}
\label{def:Wa}
\end{equation}
 \item The matrix $R(u)$ is {\em Baxter's $R$-matrix} defined by
\begin{equation}
    R(u) = \sum_{a=0}^3 W_a^R(u) \sigma^a \tensor \sigma^a,\qquad
    W_a^R(u) := W_a^L(u + \eta).
\label{def:R}
\end{equation}
 \item The indices designate the spaces on which operators act
       non-trivially. For example,
\begin{equation*}
    L_{12}(u) =
    \sum_{a=0}^3
    W_a^L(u) S^a \tensor \sigma^a \tensor 1,\qquad
    R_{23}(u) = 
    \sum_{a=0}^3
    W_a^R(u) 1 \tensor \sigma^a \tensor \sigma^a.
\end{equation*}
\end{itemize} 

Although the relation \eqref{RLL} contains parameters $u$ and $v$, the
commutation relations among $S^a$ ($a=0, \dots, 3$) do not depend on
them:
\begin{equation}
    [S^\alpha, S^0    ]_- =
    -i J_{\alpha,\beta} [S^\beta,S^\gamma]_+, \qquad
    [S^\alpha, S^\beta]_- =
                      i [S^0,    S^\gamma]_+,
\label{comm_rel}
\end{equation}
where $(\alpha, \beta, \gamma)$ stands for an arbitrary cyclic
permutation of $(1,2,3)$ and $[A,B]_\pm$ are the (anti-)commutator
$AB\pm BA$. The structure constants
$
    J_{\alpha,\beta}
    =
    ((W^L_\alpha)^2-(W^L_\beta)^2)/((W^L_\gamma)^2-(W^L_0)^2)
$
depend on $\tau$ and $\eta$ but not on $u$.

Let $l$ be a positive half integer. The {\em spin $l$ representation}
$\rho^{l}$ of the Sklyanin algebra is defined as follows: The
representation space is a space of entire functions,
\begin{multline}
    \Theta^{4l+}_{00} :=
    \{f(z) \, |\, \\
     f(z+1) = f(-z) = f(z), f(z+\tau)=\exp^{-4l\pi i(2z+\tau)}f(z) \},
\label{def:theta-space}
\end{multline}
which is of dimension $2l+1$. The generator $S^a$ of the Sklyanin
algebra acts as a difference operator on this space:
\begin{equation}
    (\rho^l(S^a) f)(z) =
    \frac{s_a(z-l\eta)f(z+\eta)-s_a(-z-l\eta)f(z-\eta)}
         {\theta_{11}(2z,\tau)},
\end{equation}
where
\begin{alignat*}{2}
    s_0(z) &=  \theta_{11}(\eta,\tau) \theta_{11}(2z,\tau),\qquad&
    s_1(z) &=  \theta_{10}(\eta,\tau) \theta_{10}(2z,\tau),\\
    s_2(z) &= i\theta_{00}(\eta,\tau) \theta_{00}(2z,\tau),\qquad&
    s_3(z) &=  \theta_{01}(\eta,\tau) \theta_{01}(2z,\tau).
\end{alignat*}
In the simplest case $l= 1/2$, $\rho^{1/2}(S^a)$ are expressed by the
Pauli matrices $\sigma^a$. We can identify $\Theta^{2+}_{00}$ and
$\Comp^2$ by
\begin{equation}
\begin{aligned}
    \theta_{00}(2z,2\tau)-\theta_{10}(2z,2\tau) 
    &\longleftrightarrow
    \begin{pmatrix} 1 \\ 0 \end{pmatrix},
\\
    \theta_{00}(2z,2\tau)+\theta_{10}(2z,2\tau) 
    &\longleftrightarrow
    \begin{pmatrix} 0 \\ 1 \end{pmatrix}.
\end{aligned}
\label{rep:identify}
\end{equation}
Under this identification $S^a$ have matrix forms
\begin{equation}
    \rho^{1/2}(S^a) = \theta_{11}(2\eta,\tau) \sigma^a.
\label{rep:pauli}
\end{equation}

The representation space $\Theta^{4l+}_{00}$ has a natural Hermitian
structure defined by the following {\em Sklyanin form}:
\begin{equation}
    \left< f(z), g(z)\right>
    :=
    \int_0^1 dx \int_0^{\tau/i} dy\,
    \overline{f(z)} g(z) \mu(z,\bar z),
\label{skl-form}
\end{equation}
where $z=x+iy$ and the kernel function $\mu(z,w)$ is defined by
\begin{equation}
    \mu(z,w)
    :=
    \frac
    {\theta_{11}(2z,\tau) \theta_{11}(2w,\tau)}
    {\displaystyle
     \prod_{j=0}^{2l+1}
     \theta_{00}(z+w+(2j-2l-1)\eta,\tau)\theta_{00}(z-w+(2j-2l-1)\eta,\tau)
    }.
\label{def:mu}
\end{equation}
The most important property of this sesquilinear positive definite 
scalar product is that the generators $S^a$ of the Sklyanin algebra
become self-adjoint:
\begin{equation}
    (S^a)^* = S^a,\text{ namely, }
    \langle f(z), S^a g(z) \rangle = \langle S^a f(z), g(z) \rangle.
\label{Sa:self-adj}
\end{equation}

In \cite{skl:83} Sklyanin also defined involutive automorphisms:
\begin{equation}
    X_a: (S^0, S^a,  S^b,  S^c) \mapsto
         (S^0, S^a, -S^b, -S^c),
\label{def:X-a}
\end{equation}
for $a=1,2,3$, where $(a,b,c)$ is a cyclic permutation of $(1,2,3)$. The
unitary operators $U_a$ defined by
\begin{equation}
\begin{aligned}
    U_1: &\Theta^{4\ell+}_{00} \owns f(z) \mapsto &
         &(U_1 f)(z) = e^{\pi i \ell} f\left(z + \frac{1}{2} \right),
\\
    U_3: &\Theta^{4\ell+}_{00} \owns f(z) \mapsto &
         &(U_3 f)(z) = e^{\pi i \ell} e^{\pi i \ell (4z+\tau)}
                       f\left(z + \frac{\tau}{2} \right),
\end{aligned}
\label{def:Ua}
\end{equation}
and $U_2= U_3 U_1$, intertwine representations $\rho^\ell \circ X_a$ and
$\rho^{\ell}$: $\rho^\ell(X_a(S^b)) = U_a^{-1} \rho^\ell(S^b)
U_a$. Operators $U_a$ satisfy the relations: $U_a^2 = (-1)^{2\ell}$,
$U_a U_b = (-1)^{2\ell} U_b U_a = U_c$.

\section{Proof of \lemref{lem:factorise-Phi}}
\label{app:factorisation}

In this appendix we prove \lemref{lem:factorise-Phi}, using the results
by Rosengren, \cite{ros:04} and \cite{ros:07}. (See also \cite{kon:05}
for the detailed proof of formulae in the elliptic case.)

Rosengren \cite{ros:07} introduced the vectors of
$\Theta^{2N+}_{00}$,
\begin{equation}
    e^N_k(z;a,b) = [z;a]_k [z;b]_{N-k},
\label{natural-basis}
\end{equation}
where $k=0,1,\dotsc,N$, $a,b\in\Comp$ and 
\begin{equation}
    [z;a] := [z+a][-z+a], \quad
    [z;a]_k := [z+a]_k [-z+a]_k.
\label{def:[z;a]}
\end{equation}
They form a basis of $\Theta^{2N+}_{00}$ when the parameters $a$ and $b$
satisfy the genericity condition:
\begin{gather*}
    a-b+2j\eta\not\in\Integer+\tau\Integer \quad
    (j=1-N,2-N,\dotsc,N-1),
\\
    a+b+2j\eta\not\in\Integer+\tau\Integer \quad
    (j=0,1,\dotsc,N-1).
\end{gather*}
Actually they are the same as the intertwining vectors
\eqref{intertwining-vec} up to parametrisation:
\begin{equation}
    \phi_{\lambda,\lambda'}(u;v)
    =
    e^{2l}_{l+m}
    \left(z;
     \frac{\lambda+ u-v}{2} + (-l+1)\eta,
     \frac{\lambda'+u-v}{2} + (-l+1)\eta
    \right).
\label{int-vec=nat-bas}
\end{equation}

Rosengren studied the change of basis $\{e^N_k(z;a,b)\}_{k=0,\dotsc,N}$
to $\{e^N_k(z;c,d)\}_{k=0,\dotsc,N}$ in \cite{ros:07} and computed the
Sklyanin form among them in \cite{ros:04}. We need the following facts
from his results.

The coefficients $R^l_k(a,b,c,d;N)$ ({\em elliptic $6j$-symbols}) of the
expansion ((5.4) in \cite{ros:07}; (3.4) in \cite{kon:05})
\begin{equation}
    e^N_k(z;a,b)
    =
    \sum_{l=0}^N
    R^l_k(a,b,c,d;N) e^N_l(z;c,d)
\label{e=Re}
\end{equation}
is expressed by the coefficients ({\em elliptic binomial coefficients})
of the expansion (the elliptic version of (3.9) in \cite{ros:07}; (3.2)
in \cite{kon:05})
\begin{equation}
    [z;a]_k
    =
    \sum_{n=0}^k C_n^k(a,b,c) [z;b]_n [z;c]_{k-n},
\label{ell-binomial}
\end{equation}
as follows (the elliptic version of the equation before Theorem 3.3 in
\cite{ros:07}; the last equation in \S4.3 of \cite{kon:05}):
\begin{equation}
    R_k^l(a,b,c,d;N)
    =
    \sum_{j=0}^{\min(k,l)}
    C^k_j(a,c,b+2(N-k)\eta) C^{N-j}_{l-j}(b,c+2j\eta,d).
\label{R=CC}
\end{equation}
The explicit expression of $C_n^k(a,b,c)$ is also known (the elliptic
version of (3.14) in \cite{ros:07}; (3.3) in \cite{kon:05}):
\begin{multline}
    C_n^k(a,b,c)
\\
    =
    \frac{[2\eta]_k}
         {[2\eta]_n [2\eta]_{n-k}}
    \frac{[a-c]_n [a+c+2(k-n)\eta]_n [a-b]_{k-n} [a+b+2n\eta]_{k-n}}
         {[b-c+2(n-k)\eta]_n [c-b-2n\eta]_{k-n} [b+c]_k}.
\label{C}
\end{multline}
What is necessary for us later is the extremal case
$R_N^N(a,b,c,d;N)$. Using \eqref{R=CC} and \eqref{C}, we have
\begin{equation*}
 \begin{split}
    R_N^N(a,b,c,d;N)
    &=
    \sum_{j=0}^{N}
    C^N_j(a,c,b) C^{N-j}_{N-j}(b,c+2j\eta,d)
\\
    &=
    \sum_{j=0}^{N}
    C^N_j(a,c,b)
    \frac{[b-d]_{N-j} [b+d]_{N-j}}
         {[c+2j\eta-d]_{N-j} [c+2j\eta+d]_{N-j}}
\\
    &=
    \frac{1}{[c-d]_N [c+d]_N} \times
\\
    &\times
    \sum_{j=0}^{N}
    C^N_j(a,c,b)
    [c-d]_j [c+d]_j [b-d]_{N-j} [b+d]_{N-j}.
 \end{split}
\end{equation*}
Comparing the last expression with \eqref{ell-binomial}, we obtain the
formula:
\begin{equation}
    R_N^N(a,b,c,d;N) = \frac{[d;a]_N}{[d;c]_N}.
\label{RNN}
\end{equation}
(It is natural that $b$ does not appear in the right hand side, since
$e^N_N(z;a,b)$ does not depend on $b$.)

We also need the following orthogonality relation in
\cite{ros:04}\footnote{$\langle f, g\rangle$ in \cite{ros:04} is
$\langle g, f\rangle$ in \cite{skl:83}. It is linear in $f$ and
conjugate linear in $g$. We follow the convention in \cite{skl:83}.}
(Theorem 3.4).
\begin{equation}
 \begin{split}
    &\langle 
     e^N_l(z;-\bar d + (1-N)\eta + \tfrac{\tau+1}{2},
             -\bar c + (1-N)\eta - \tfrac{\tau+1}{2}),
     e^N_k(z;c,d)
    \rangle
\\
    &=
    C_N e^{2\pi i(-dl + c(N-l) - N(1+\tau)/4)}
    \Gamma^N_k(c,d) \delta_{k,l},
 \end{split}
\label{<el,ek>}
\end{equation}
where
\begin{equation*}
    C_N
    =
    \frac{-2\eta e^{3\pi i \tau/4}}
         {[2(N+1)\eta] \prod_{j=1}^\infty (1-e^{2j \pi i \tau})^3}
\end{equation*}
is a constant depending only on $(\tau,\eta,N)$ and
\begin{equation}
 \begin{split}
    \Gamma^N_k(c,d)
    &=
    e^{\pi i N(\tau-1)/2}
    \frac{[c-d-2N\eta]}{[c-d+2(2k-N)\eta]} \times
\\
    &\times
    \frac{[2\eta]_k [c-d+2\eta]_k}
         {[-2N\eta]_k [c-d-2N\eta]_k}
    [c-d+2(1-N)\eta]_N [c+d]_N.
 \end{split}
\end{equation}
In particular,
\begin{equation}
 \begin{split}
    \Gamma^N_N(c,d)
    &=
    e^{\pi i N(\tau-1)/2}
    \frac{[c-d-2N\eta]}{[c-d+2N\eta]} \times
\\
    &\times
    \frac{[2\eta]_N [c-d+2\eta]_N}
         {[-2N\eta]_N [c-d-2N\eta]_N}
    [c-d+2(1-N)\eta]_N [c+d]_N
\\
    &=
    e^{\pi i N(\tau+1)/2} [c-d]_N [c+d]_N.
 \end{split}
\label{GammaNN}
\end{equation}
Using these formulae, let us compute the Sklyanin form
$
    \langle
     e^N_N(z;\alpha,\beta), e^N_N(z;\gamma,\delta)
    \rangle
$.
In the expansion of the form \eqref{e=Re},
\begin{equation}
    e^N_N(z;\gamma,\delta)
    =
    \sum_{n=0}^N R^n_N(\gamma,\delta,c,d;N)
    e^N_n(z;c,d),
\label{e(gamma,delta)=Re(c,d)}
\end{equation}
of the second factor, we choose the parameters $c$ and $d$, so that
$\{e^N_n(z;c,d)\}_{n=0,\dotsc,N}$ is a dual basis (up to normalisation)
to $\{e^N_n(z;\alpha,\beta)\}_{n=0,\dotsc,N}$. According to
\eqref{<el,ek>},
\begin{equation*}
    \alpha =-\bar d +(1-N)\eta + \frac{\tau+1}{2}, \quad
    \beta  =-\bar c +(1-N)\eta - \frac{\tau+1}{2},
\end{equation*}
namely,
\begin{equation}
    c = -\bar \beta  +(1-N)\eta - \frac{-\tau+1}{2}, \quad
    d = -\bar \alpha +(1-N)\eta + \frac{-\tau+1}{2}.
\label{c,d}
\end{equation}
Applying $\langle e^N_N(z;\alpha,\beta), \cdot\rangle$ to
\eqref{e(gamma,delta)=Re(c,d)} and using the orthogonality relation
\eqref{<el,ek>}, we have
\begin{equation*}
 \begin{split}
    &\langle
     e^N_N(z;\alpha,\beta), e^N_N(z;\gamma,\delta)
    \rangle
\\
    &=
    C_N e^{2\pi i (-d N - N(1+\tau)/4)} \Gamma^N_N(c,d)
    R^N_N(\gamma,\delta,c,d;N)
\\
    &=
    C_N e^{2\pi i (-d N - 2\pi i N(1+\tau)/4)}
        \bigl(e^{\pi i N(\tau+1)/2} [c-d]_N [c+d]_N\bigr)
    \left(\frac{[d;\gamma]_N}{[d;c]_N}\right)
\\
    &=
    C_N e^{-2\pi i dN} [d;\gamma]_N.
 \end{split}
\end{equation*}
Substituting \eqref{c,d},
\begin{equation}
 \begin{split}
    &\langle
     e^{N}_{N}(z;\alpha,\beta), e^{N}_{N}(z;\gamma,\delta)
    \rangle
\\
    ={}&
    C_{N} e^{-2\pi i N (-\bar\alpha + (1-N)\eta + (-\tau+1)/2)}
\\
    \times&\prod_{j=0}^{N-1}
    [-\bar\alpha + (1-N)\eta + \tfrac{-\tau+1}{2} + \gamma+2j\eta]
    [ \bar\alpha - (1-N)\eta - \tfrac{-\tau+1}{2} + \gamma+2j\eta]
\\
    ={}&
    C_{N} e^{\pi i N \tau/2}
    \prod_{j=0}^{N-1}
    \theta_{00}(\gamma-\bar\alpha + (2j-N+1)\eta,\tau)
    \theta_{00}(\gamma+\bar\alpha + (2j+N-1)\eta,\tau).
 \end{split}
\label{<e,e>}
\end{equation}

\medskip
In order to prove \lemref{lem:factorise-Phi}, we need to compute the
Sklyanin form of two local pseudo vacuum vectors,
$
    \omega_{\sigma\lambda}(-\bar u;\sigma v)
    =
    \phi_{\sigma\lambda,\sigma\lambda+4l\eta}(-\bar u;\sigma v)
$
and
$
    \omega_{\sigma'\lambda'}(u';\sigma' v')
    =
    \phi_{\sigma'\lambda',\sigma'\lambda'+4l\eta}(u';\sigma' v')
$.
We identify them with the vectors $e^{2l}_{2l}(z;a,b)$ by the formula
\eqref{int-vec=nat-bas}, but not directly. For the later purpose, we
should first modify the expression of $\omega_{\sigma\lambda}(u;v)$,
\begin{equation*}
    \omega_{\sigma\lambda}(u;\sigma v)
    =
    \prod_{j=0}^{2l-1}
    [ z+\tfrac{\sigma\lambda + u - \sigma v}{2} + (2j-l+1)\eta]
    [-z+\tfrac{\sigma\lambda + u - \sigma v}{2} + (2j-l+1)\eta].
\end{equation*}
Since $[z]$ is an odd function and $\sigma=\pm1$, it is rewritten as
\begin{equation*}
    \omega_{\sigma\lambda}(u;\sigma v)
    =
    \prod_{j=0}^{2l-1}
    [ z+\tfrac{\lambda + \sigma u - v}{2} + \sigma (2j-l+1)\eta]
    [-z+\tfrac{\lambda + \sigma u - v}{2} + \sigma (2j-l+1)\eta].
\end{equation*}
The set $\{\sigma(2j-l+1)\eta\mid j=0,\dotsc,2l-1\}$ is equal to
$\{\sigma l \eta -(2l-1)\eta + 2j\eta\mid j=0,\dotsc,2l-1\}$. Therefore
we can identify $\omega_{\sigma\lambda}(u;v)$ with $e^{2l}_{2l}$ as
follows: 
\begin{equation}
    \omega_{\sigma\lambda}(u;\sigma v)
    =
    e^{2l}_{2l}(z;
    \tfrac{\lambda+\sigma u - v}{2} + \sigma l\eta -(2l-1)\eta, \ast
    ),
\label{omega=e}
\end{equation}
Here the second parameter in $e^{2l}_{2l}$ is irrelevant.

\medskip
The final step of the proof is the computation of
$    \langle
    \omega_{\sigma\lambda}(-\bar u;\sigma v),
    \omega_{\sigma'\lambda'}(u';\sigma' v')
    \rangle.
$
Taking \eqref{omega=e} into account, we substitute $N=2l$ and
\begin{align}
    \alpha
    &=
    \tfrac{\lambda-\sigma \bar u  - v }{2} + \sigma  l\eta -(2l-1)\eta,
\label{alpha}
\\
    \gamma
    &=
    \tfrac{\lambda'+\sigma'    u' - v'}{2} + \sigma' l\eta -(2l-1)\eta,
\label{gamma}
\end{align}
into \eqref{<e,e>}, which gives
\begin{equation*}
 \begin{split}
    &\langle
    \omega_{\sigma\lambda}(-\bar u;\sigma v),
    \omega_{\sigma'\lambda'}(u';\sigma' v')
    \rangle
\\
    ={}&
    C_{2l} e^{\pi i l \tau}
    \prod_{j=0}^{2l-1}
    \theta_{00}(
      \tfrac{\lambda'-\bar\lambda}{2} 
    + \tfrac{\sigma'u' + \sigma u}{2} 
    + (\sigma'-\sigma) l\eta
    + \tfrac{-v' + \bar v}{2} 
    + (2j-N+1)\eta,
    \tau)
\\
    \times&
    \prod_{j=0}^{2l-1}
    \theta_{00}(
      \tfrac{\lambda'+\bar\lambda}{2} 
    + \tfrac{\sigma'u' - \sigma u}{2} 
    + (\sigma'+\sigma) l\eta - 2(2l-1)\eta
    + \tfrac{-v' - \bar v}{2} 
    + (2j-N+1)\eta,
    \tau).
 \end{split}
\end{equation*}
This is one of the factors in \eqref{Phi=prod} and has the desired
factorised structure as in \eqref{Phi=FG}.
\qed

\end{appendices}


\end{document}